\documentclass[12pt]{article}

\usepackage{amsfonts,epsfig}

\newcommand{\ptl}{\partial}
\newtheorem{proposition}{Proposition}
\newtheorem{theorem}{Theorem}

\begin{document}
\title{Wiener--Hopf matrix factorization  using ordinary differential equations
in the commutative case}

\author{A.V. Shanin}
\date{}

\maketitle

\begin{abstract}
A matrix factorization problem is considered. The matrix to be factorized
is algebraic, has dimension $2\times 2$ and belongs to  Moiseev's class.
A new method of factorization is proposed. First, the matrix factorization problem
is reduced to a Riemann--Hilbert problem using the Hurd's method. Secondly,
the Riemann--Hilbert problem is embedded
into a family of Riemann--Hilbert problems indexed by a variable $b$ taking values on a half--line.
A linear ordinary differential equation (ODE1) with respect to $b$ is derived. The coefficient
of this equation remains unknown at this step.
Finally, the coefficient of the ODE1 is computed. For this, it is proven that it obeys
a non-linear ordinary differential equation (ODE2) on a half--line. Thus, the numerical procedure
of matrix factorization becomes reduced to two runs of solving of ordinary differential equations
on a half--line: first ODE2 for the coefficient of ODE1, and then ODE1 for the unknown function.
The efficiency of the new method is demonstrated on some examples.
\end{abstract}

\section{Introduction}

Many diffraction problems can be  transformed into matrix factorization
problems \cite{Noble}. Typically, these diffraction problems are
2D problems with different boundaries occupying positive
and negative parts of the $x$-axis. There emerges a known matrix $G$ analytical in a thin strip
going along the real axis of a complex variable $k$, and it is necessary to represent it as a product
\begin{equation}
G(k) = U^{-1} (k) W(k) ,
\label{eq0101}
\end{equation}
where $W$ and $U$ are matrices analytical and having no zeros of the determinant in the upper
and lower half--plane, respectively. Also, both matrices should have algebraic growth in corresponding half--planes.

In the scalar case, which can be considered as a degenerate case of $1\times 1$ matrix, the solution can
be readily achieved by taking the logarithm of the matrix and performing the additive decomposition by means of
Cauchy's integral and Sokhotsky's formula \cite{Noble}.
Returning to the matrices of order $N > 1$,
this approach can be generalized for Moiseev's matrices \cite{Moiseev} having form
\begin{equation}
G = \sum_{n=0}^{N-1} g_n(k) A^n (k),
\label{eq0102}
\end{equation}
where $g_n(k)$ are scalar function, and $A$ is a polynomial matrix.
In the simplest case of matrix $A$ having distinct eigenvalues almost everywhere,
$A$ can be decomposed as
\begin{equation}
A = T D T^{-1},
\label{eq0103}
\end{equation}
where $T$ is the matrix of the eigenvectors and $D$ is a diagonal matrix composed of the
eigenvalues. Both $T$ and $D$ are algebraic matrices, and one can introduce the Riemann surface
$\cal R$ on which $T$ and $D$ are single--valued. Further, the matrix factorization problem
becomes reduced to a scalar Riemann--Hilbert problem on~$\cal R$. This problem can be solved in terms of Abelian
integrals with the help of Jacobi's inversion problem \cite{Zverovich}.
So, the solution of the problem of factorization of (\ref{eq0102}) is known at least formally,
and it possibly can be used for practical needs. Some examples can be found e.g.\ in \cite{Antipov}.
Simpler, but more popular cases \cite{Khrapkov,Jones} can be described
as particular cases of (\ref{eq0102}).
Khrapkov's method \cite{Khrapkov} is rather simple and leads to straightforward computations, but
for a broad class of matrices it produces
non-algebraic growth at infinity. Moiseev's method can be considered as
a remedy enabling one to avoid this growth. Another technique to avoid the non-algebraic growth has been
proposed in \cite{Daniele}. This technique also includes some numerical stages.
A review of the commutative factorization and a development of ideas of \cite{Jones} can be found
in~\cite{AbrahamsCom}.

If $G$ cannot be represented as (\ref{eq0102})
then some numerical \cite{Gakhov} or approximate (e.g.\ \cite{Abrahams}) methods can be applied.

In the current work we consider matrices of Moiseev's class (\ref{eq0102}) and develop a new technique
which is arguably simpler in practical realization than the Moiseev--Zverovich or Daniele procedure.
The new technique can be applied only when the $g_n(k)$ in (\ref{eq0102}) are algebraic functions.
This is an important restriction, however in much of the practical situation this restriction is
fulfilled.
The new procedure comprises three steps.
Firstly, the matrix factorization problem
is reduced to a Riemann--Hilbert problem using the Hurd's method \cite{Hurd}.
Namely, instead of studying the factors $W$ and $U$ we are studying only the factor
$U$ continued into the upper half--plane of~$k$. The $k$--plane is cut along half--lines connecting
the branch points of $G$ located in the upper half--plane, namely the points $k_j$,
with $k = + i \infty$.
A Hilbert problem is formulated on the half--lines $(k_j, k_j + i \infty)$.
As Hurd mentioned, the new problem can be simpler than the initial matrix factorization problem.

At the second step the Riemann--Hilbert problem is embedded
into a family of Riemann--Hilbert problems indexed by a variable $b$ taking values on a half--line.
Namely, for the whole family the coefficients $H_j(k)$ remain the same, but the contours on which the
functional equations should be fulfilled are changed from $(k_j , k_j + i \infty)$ to
$(k_j + b, k_j + i \infty)$, where $b$ is an imaginary number taking values from $0$ to $+i \infty$.
Thus, we can define the family of solutions $U(b, k)$. The solution of the initial problem is denoted by $U(0,k)$.
A linear ordinary differential equation (ODE1) with respect to $b$ is derived for $U(b,k)$.
The coefficient  of this equation remains unknown on this step.

Finally, step the coefficient of the ODE1 is computed. For this, it is proved that it obeys
a non-linear ordinary differential equation (ODE2) on a half--line. Thus, the numerical procedure
of matrix factorization becomes reduced to  solving two ordinary differential equations
on a half--line: first ODE2 for the coefficient of ODE1, and then ODE1 for the unknown function.

Some numerical results are presented. Namely, we demonstrate that the new procedure applied to
a
matrix belonging to the Khrapkov's class is factorized exactly the same way as by the traditional Khrapkov's procedure. Moreover, we apply our method to the matrix emerging in \cite{Antipov}.

\section{Problem formulation and Hurd's procedure}

Let $G(k)$ be an algebraic matrix $N\times N$ having no singularities and no zeros of determinant on the
real axis and tending to the unit matrix $I$ of dimension $N\times N$ as $|k| \to \infty$.
 Our aim is to find the decomposition (\ref{eq0101}) valid in some strip
 $|{\rm Im}[k]| < \epsilon$ with $U$  having no singularities or zeros of the determinant in the lower half--plane and on the real axis,  and $W$ having no singularities or zeros of the determinant in the upper half--plane, maybe except several points, where
 poles or zeros are allowed.
We demand that the unknown functions $U$ and $W$ tend to $I$ as $|k| \to \infty$.
Some restrictions on $G$ will be imposed below.

\begin{figure}[ht]
\centerline{\epsfig{file=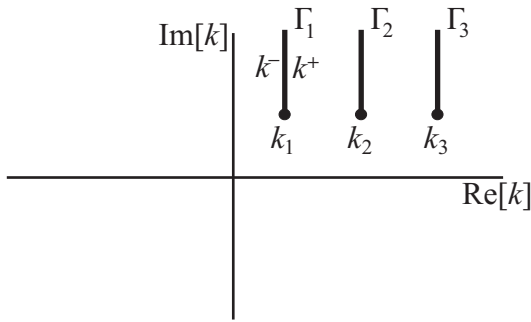}}
\caption{Cuts $\Gamma_j$} \label{fig01}
\end{figure}

Apply Hurd's procedure \cite{Hurd} as follows.
Let $k_j$, $j = 1,\dots,p$ be branch points of matrix $G$ in the upper half--plane. Connect the points
$k_j$ with $i\infty$ by the cuts
\[
\Gamma_j = (k_j , k_j + i \infty)
\]
parallel to the imaginary axis.
Let contours $\Gamma_j$ do not pass through other branch points, poles or zeros of the
determinant of~$G$.
Continue function $U(k)$ into the upper half--plane cut along the lines $\Gamma_j$
by the relation
\begin{equation}
U(k) \equiv W(k) G^{-1}(k).
\label{eq0201}
\end{equation}
Note that $W$ is defined and regular in the upper half--plane and $G$ is defined in the
upper half--plane with the cuts~$\Gamma_j$.
Define by $U(k^+)$ and $U(k^-)$ for $k \in \Gamma_j$ the values of $U$ taken on the right and on the
left shore of $\Gamma_j$, respectively (see Fig.~\ref{fig01}). Similarly, define the values $G(k^+)$
and $G(k^-)$. For some $k \in \Gamma_j$
\[
U(k^+) \equiv W(k) G^{-1}(k^+),
\]
\[
U(k^-) \equiv W(k) G^{-1}(k^-)
\]
(note that $W(k^-)= W(k^+) = W(k)$).
Then,
\begin{equation}
U(k^+) = U(k^-) H_j (k),
\label{eq0202a}
\end{equation}
\begin{equation}
 H_j(k)\equiv G(k^-) G^{-1}(k^+), \qquad k \in \Gamma_j
\label{eq0202}
\end{equation}
The set of equations (\ref{eq0202}) taken for $j = 1,\dots,p$ constitute the Riemann--Hilbert
problem in Hurd'd formulation. It was Hurd's observation that this problem can be simpler than the initial matrix factorization problem.

In our case each coefficient $H_j(k)$ of the Riemann--Hilbert problem can be continued analytically into some half--strip $\Omega + k_j$
\[
\Omega = \{ {\rm Re}[k]< \epsilon, \, {\rm Im}[k] > 0 \}.
\]
Note that the point $k_j$ does not belong to $\Omega + k_j$.
Also $H_j(k)$ are algebraic functions, so each of them can be continued onto some Riemann surface.

Let all $H_j(k)$ tend to $I$  ${\rm Im}[k] \to \infty$
(this restriction is fulfilled if $G(k)\to I$ on all sheets).
We forget about $W(k)$ and look for $U(k)$ on the complex plane $k$
cut along the contours $\Gamma_j$ having no singularities and no zeros of the determinant on the cut plane, obeying the problem (\ref{eq0202a}). We recall that $U(k) \to I$
as $|k| \to \infty$.

The behavior of $U(k)$ at the points $k_j$ will be specified below in such a
way that the problem possesses a unique solution.

We assume also that matrices $H_j$ have distinct eigenvalues almost everywhere.

All restrictions described above correspond to a quite general matrix factorization problem and are
easy to fulfil by, e.g.\ slight change of the contour position. If matrices $H_j$ do not tend to $I$ the method can be easily modified also. Here we are going to pose the most strong restriction: we assume that
all branches of all matrices $H_j(k)$ taken for arbitrary {\em affix}\/ $k$ commute with each other, i.e.
for each $k$
\begin{equation}
H_{j_1} (k) H_{j_2} (k) = H_{j_2} (k)  H_{j_1} (k),
\label{eq0203}
\end{equation}
where $H_{j_1}(k)$ and $H_{j_2}(k)$ are any possible continuations of $H_{j_1}$ and $H_{j_2}$ to $k$.
The meaning of this restriction is discussed in the next section.

\section{Functional--commutative and branch--com\-mu\-ta\-ti\-ve matrices}

All existing analytical
approaches to matrix factorization are available only for matrices admitting a commutative
factorization, i.e.\ a representation of the form
\begin{equation}
G(k) = U^{-1} (k) W(k) = W(k) U^{-1}(k).
\label{eq0301}
\end{equation}
The theory of commutative
matrix factorization starts from \cite{Chebotarev} where a concept of functional--commutative matrix has been introduced. A functional--commutative matrix is a matrix commuting with its singular integral. The property
of functional--commutativity is not easy to check for an arbitrary matrix. That is why, in \cite{arxiv1} we
introduced branch--commutative matrices. Namely, an algebraic matrix $G(k)$ is called branch-commutative if
for any $k$ the values $G_j(k)$ corresponding to different branches of $G$ commute with each other.

To formulate the main result of \cite{arxiv1} we need one more definition. A Riemann surface of an algebraic matrix is called
{\em balanced}\/ if each sheet of it can be reached from any fixed sheet only by bypassing the
branch points located in the upper half--plane and only bypassing the branch points lying in the negative
half--plane. Most of the known matrices arising in practical problems have balanced Riemann surfaces.

The main result of \cite{arxiv1} is as follows. If an algebraic matrix $G$ with balanced Riemann surface
admits commutative factorization then it is branch--commutative. Vice versa, a branch--commutative matrix
can be represented in the form (\ref{eq0102}), and thus the Moiseev's method can be applied to it.

Note that if matrix $G$ is branch--commutative then the property of (\ref{eq0203}) for the matrices $H_j$
defined by (\ref{eq0202}) is valid. Thus, the matrices to which the method described here can be applied are
(with some unimportant restrictions) the same as the matrices, to which the Moiseev's method is applicable.

Let us formulate one important consequence of the property (\ref{eq0203}).

\begin{proposition}
If property (\ref{eq0203}) is fulfilled then there exists rational matrix $B(k)$ commuting with all
matrices $H_j (k)$.
\label{proposition1}
\end{proposition}

The proof of the proposition is as follows. Represent $H_1$ in the form
\begin{equation}
H_1 (k) = P(k) F_1(k) P^{-1}(k),
\label{eq0302}
\end{equation}
where $P(k)$ is the matrix composed of the eigenvectors of $H_1$ normalized, say, by making the first
component of each vector equal to~1. Respectively, $F_1(k)$ is a diagonal matrix composed of
scalar functions $f_1(k), \dots , f_N (k)$.

It is known that if two matrices commute then normalized eigenvectors of the matrices coincide \cite{Lancaster}.
Matrix $P$ is algebraic, so it is single--valued on some Riemann surface. Since the values of
$H_j(k)$ taken on different sheets (with the same $k$) commute, we can conclude that when a branch point of $P$
is bypassed the columns of $P$ are just permuted.

Construct matrix $B$ in the form
\begin{equation}
B(k) = P(k) D(k) P^{-1}(k),
\label{eq0303}
\end{equation}
where $D(k)$ is a diagonal matrix with the scalar functions $h_1(k), \dots , h_N(k)$ on the diagonal.
Let the functions $h_m$ be branches of some algebraic function $h$ having the same branch
points as $P$. Moreover, let the values $h_m$ be permuted the same way as the columns of $P$ when
the branch points are bypassed. Then the function $B(k)$ is single--valued, and therefore rational.

A proper choice of the functions $h_m$ is as follows:
\begin{equation}
h_m = \sum_{n = 1}^N \beta_n(k) P_{n,m} (k),
\label{eq0304}
\end{equation}
where $\beta_n(k)$ is an arbitrary set of scalar rational functions (provided none of $h_m$ is identically
zero). Values of $P_{n,m}$ are elements of matrix~$P$.

Since $B$ commutes with $H_1$ and all matrices $H_j$ have distinct eigenvalues almost everywhere,
matrix $B$ commutes with every $H_j$.
Note that all other matrices $H_j$ can be represented as
\begin{equation}
H_j (k) = P(k) F_j(k) P^{-1}(k).
\label{eq0302a}
\end{equation}

Due to arbitrariness of the choice of $\beta_n (k)$ one can make matrix $B$ having simple poles
and tending to $I$ as $|k| \to \infty$.
The matrix $B$ possessing all these properties plays an important role below.

\section{A family of Riemann--Hilbert problems and derivation of ODE1}

\subsection{Family of Riemann--Hilbert problems}

We have reduced the matrix factorization problem to finding the
function $U(k)$ obeying equations (\ref{eq0202a}) on the cuts $\Gamma_j$.
To solve this problem we use the idea described in detail in
\cite{arxiv2}. Namely, we are fixing the functions $H_j(k)$, defined
and continuous on contours $\Gamma_j$ (and regular in the strips $\Omega + k_j$),
and introduce truncated contours
\[
\Gamma_j(b) = (k_j + b , k_j + i \infty),
\]
where $b$ is an imaginary number $b \in (0, i \infty)$. Consider a family of problems
for the function $U(b, k)$
set by the relations
\begin{equation}
U(b,k^+) = U(b,k^-) H_j (k), \qquad k \in \Gamma_j(b).
\label{eq0401}
\end{equation}
We assume that for each $b$ the matrix function $U(b,k)$ is single--valued,
continuous, and free of zeros of determinant on the plane of $k$ cut along the contours $\Gamma_j (b)$.
It
tends to $I$ as $|k| \to \infty$. We assume also that the behavior of $U(b,k)$
at the points $k_j + b$ is derived by continuity from the conditions
formulated for large ${\rm Im}[b]$ (see below). Obviously,
\begin{equation}
U(k) = U(0, k).
\label{eq0402}
\end{equation}
The main idea of the method is to study the behavior of $U(b, k)$ as a function of $b$.

Embedding of $U(k)$ into a family $U(b,k)$ enables us to define behavior of $U(k)$ at
$k = k_j$ in the most natural way. Expand matrices $H_j$ in the form (\ref{eq0302a}).
The leading term of $U(b,k)$ near the point $k =  k_j + b$ can be written in the form
\begin{equation}
U(b,k) \approx K_j(k_j + b) (k - (k_j + b))^{\log(F_j(k_j+ b)) / (2 \pi i)} K_j^{-1}(k_j + b)
\label{eq0408}
\end{equation}
where $F_j$ are taken from (\ref{eq0302}), and $K_j$ are some non-singular matrices.
The branch of the logarithm should be fixed as follows. For $b \to i\infty$
the matrices $F_j(k_j + b)$ tend to~$I$. For these values we choose the branch of
logarithm close to zero matrix. Then, for other values of $b$ we choose the branch of
logarithm by continuity. Such a choice enables us to avoid discussing partial indices of the
initial Riemann--Hilbert problem.

\subsection{Form of ODE1 for a single cut}

Let the number of cuts $p$ be equal to 1, i.e. let there exists only one cut $\Gamma_j$. This corresponds
to a function $G(k)$ having a single branch point in the positive half--plane.
This case has been studied in \cite{arxiv2}. Here we formulate the main theorem of \cite{arxiv2}
with a short proof.

For this, it is necessary to introduce a notation of the {\em ordered exponential}\/ (the
term comes from quantum mechanics).
Namely, let $\gamma$ be a contour (directed one) connecting the points
$\tau_1$ and $\tau_2$ ($\tau_1$ is a starting point), and let $C(\tau)$ be a $N\times N$ matrix defined
on $\gamma$. Consider a matrix equation
\begin{equation}
\frac{d}{d \tau} X (\tau) = C(\tau) X(\tau)
\label{eq0401a}
\end{equation}
taken with the initial condition $X(\tau_1) = I$.
Solve this equation along contour $\gamma$ and define the value $X(\tau_2)$.
By definition,
\begin{equation}
{\rm OE}_{\gamma} [C(\tau)d\tau] \equiv X(\tau_2).
\label{eq0401b}
\end{equation}
This notation is just a convenient way to refer to a solution of an ordinary differential equation.

\begin{figure}[ht]
\centerline{\epsfig{file=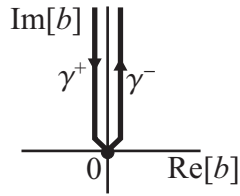}}
\caption{Contours $\gamma^+$ and $\gamma^-$}
\label{fig02}
\end{figure}

\begin{theorem}
{\bf a)} There exists $N\times N$ matrix $s_1 (b)$ analytical in the strip
$\Omega$,
such that $U(b, k)$ obeys an ordinary differential equation
(ODE1)
\begin{equation}
\frac{\ptl}{\ptl b} U(b,k) = \frac{s_1 (b)}{k - (k_1 + b)} U(b, k).
\label{eq0403}
\end{equation}
The initial condition for this equation is as follows:
\begin{equation}
U(i \infty , k) = I.
\label{eq0404}
\end{equation}

{\bf b)}
Let there exist a $N\times N$ matrix $s_1 (\tau)$ analytic in $\Omega$ and such that
\begin{equation}
{\rm OE}_{\gamma} \left[\frac{s_1(\tau)}{k - (\tau + k_1)} d \tau   \right]
= H_1 (k)
\label{eq0405}
\end{equation}
for $k \in (k_1 , k_1 + i \infty)$. Contour $\gamma$ is a concatenation of $\gamma^+$
and $\gamma^-$ (see Fig.~\ref{fig02}).
Then solution $U(b , k)$ is given by the formula
\begin{equation}
U(b , k) = {\rm OE}_{\gamma_b} \left[  \frac{s_1(\tau)}{k - (\tau + k_1)} d\tau  \right]
\label{eq0406}
\end{equation}
where contour $\gamma_b$ goes from $ i \infty$ to $b$ along $\Gamma_1$.
\end{theorem}

Let us outline the proof. Consider part a). Consider the function
\begin{equation}
S(b, k) =   \frac{\ptl U(b,k)}{\ptl b}   U^{-1} (b, k).
\label{eq0407}
\end{equation}
This function is analytic in the plane cut along $(k_j + b , k_j + i \infty)$. Consider its
behavior at the cut. Note that the coefficient $H_1$ does not depend on $b$.
Thus,
\[
\frac{\ptl U(b,k^+)}{\ptl b}  = \frac{\ptl U(b,k^-)}{\ptl b}  H_1 (k),
\quad
k \in \Gamma_1 + b.
\]
Using this relation with (\ref{eq0401}) we conclude that $S(b, k^+) = S(b, k^-)$, and therefore the
function is single--valued. According to the condition at infinity for $U$, $S$ should decay
 as $|k|\to \infty$.
The only singularity of $S$ in the finite part of the $k$--plane is $k = k_1 + b$.
The leading term of the singularity is given by (\ref{eq0408}). According to this,
$S$ has a simple pole at $k = k_1 + b$, and
\begin{equation}
S(b, k) = - \frac{1}{2 \pi i (k - (k_1 + b))}K_1(k_1+ b) \, \log(F_1(k_1+ b)) \, K_1^{-1}(k_1+ b) .
\label{eq0409}
\end{equation}
Finally,
\begin{equation}
s_1(b) = - \frac{1}{2 \pi i} K_1(k_1+ b) \, \log(F_1(k_1+ b)) \, K_1^{-1}(k_1+ b) .
\label{eq0410}
\end{equation}
The choice of branch of the logarithm has been discussed in the previous subsection.

Analyticity of the coefficient $s_1$ in $\Omega + k_1$ follows from the fact that the
cut on which functional equation (\ref{eq0401}) is set can be deformed (without changing its starting point)
arbitrarily within $\Omega + k_1$, and the solution remains the same while the contour changes.

The initial condition (\ref{eq0404}) follows from general properties of the Riemann--Hilbert
problem \cite{Gakhov}.

Consider part b) of the theorem. Let us show that (\ref{eq0406}) is a solution of the problem
(\ref{eq0401}). Note that due to analyticity of $s_1$ the contour $\gamma$ in
(\ref{eq0405}) can be deformed provided that
it does not cross the singularity $\tau = k-k_1$ of the coefficient. Let be $k \in \Gamma_1 (b)$, and thus
${\rm Im}[k] > {\rm Im}[k_1 + b]$. Deform contour $\gamma$ into $\gamma + b$:
\begin{equation}
{\rm OE}_{\gamma + b} \left[\frac{s_1(\tau)}{k - (\tau + k_1)} d \tau   \right]
= H_1 (k)
\label{eq0411}
\end{equation}
According to general properties of the ordinary differential equations and the ordered
exponential notations \cite{arxiv2},
\[
{\rm OE}_{\gamma + b} \left[\frac{s_1(\tau)}{k - (\tau + k_1)} d \tau   \right] =
\]

\begin{equation}
\left(
{\rm OE}_{b + \gamma^- } \left[\frac{s_1(\tau)}{k - (\tau + k_1)} d \tau   \right]
\right)^{-1}
{\rm OE}_{b + \gamma^+ } \left[\frac{s_1(\tau)}{k - (\tau + k_1)} d \tau   \right].
\label{eq0412}
\end{equation}
Note that according to (\ref{eq0406})
\[
{\rm OE}_{b + \gamma^+ } \left[\frac{s_1(\tau)}{k - (\tau + k_1)} d \tau   \right]=
U(b, k^+),
\]
\[
{\rm OE}_{b + \gamma^- } \left[\frac{s_1(\tau)}{k - (\tau + k_1)} d \tau   \right]=
U(b, k^-),
\]
Thus, (\ref{eq0412}) is equivalent to (\ref{eq0401}).

\subsection{Form of ODE1 for several cuts}

If there are $p > 1$ branch points $k_j$ (and, thus, several cuts $\Gamma_j$) Theorem~1 can be modified,
while the reasoning remains basically the same.
Coefficient $S$ from (\ref{eq0407}) can be proven to be single--valued and decaying, but it should
have $p$ simple poles $k = k_j + b$. Therefore equation (\ref{eq0403}) has form
\begin{equation}
\frac{\ptl}{\ptl b} U(b,k) =
\left( \sum_{j=1}^p \frac{s_j (b)}{k - (k_j + b)} \right)U(b, k).
\label{eq0501}
\end{equation}
with $p$ unknown matrices $s_j(b)$ analytical in $\Omega$. The initial conditions are the same as
for one cut (i.e. (\ref{eq0404})). The form of the solution follows from (\ref{eq0501}) and
(\ref{eq0404}):
\begin{equation}
U(b , k) = {\rm OE}_{\gamma_b} \left[
\sum_{j = 1}^p \frac{s_j(\tau)}{k - (\tau + k_j)} d\tau  \right].
\label{eq0502}
\end{equation}
A generalization of (\ref{eq0410}) has form
\begin{equation}
s_j(b) = - \frac{1}{2 \pi i} K_j(k_j+ b) \, \log[F_j(k_j+ b)] \, K_j^{-1}(k_j+ b) .
\label{eq0502a}
\end{equation}
for some unknown matrices $K_j$ and known (up to transmutations) diagonal matrices $F_j$.

Finally, condition (\ref{eq0405}) should be rewritten as
\begin{equation}
{\rm OE}_{\gamma} \left[\sum_{j = 1}^p \frac{s_j(\tau)}{k - (\tau + k_j)} d \tau   \right]
= H_m (k) , \quad k \in \Gamma_m + b, \quad m = 1 \dots p.
\label{eq0503}
\end{equation}

\section{ODE2}

\subsection{Derivation of ODE2}

Formula (\ref{eq0406}) (or (\ref{eq0502})) cannot be used immediately to find solution
$U(k) = U(0,k)$ since matrix functions $s_j(b)$ are unknown. Thus, before finding $U$ one should
find $s_j$ somehow.
In \cite{arxiv2} it has been proposed to use equation (\ref{eq0405}) to find $s_1$. A numerical
procedure has been proposed and tested.
Application of this procedure does not require branch--commutativeness, so the method is potentially applicable to a much wider class of problem than Moiseev's class.
However, this procedure is rather sophisticated and it does
not reveal the mathematical nature of the solution. Here we are proposing another technique reducing
the determination
of $s_j$ to solving a (nonlinear) ordinary differential equation. Unfortunately, the
new technique is applicable only to Riemann--Hilbert problems obeying relations (\ref{eq0203}).

The key idea of the new method is to use matrix $B(k)$ defined by (\ref{eq0303}). Namely,
we construct a rational matrix $B(k)$ commuting with $H_j(k)$ on all their sheets, behaving as
$B(k) \to I$ as $|k| \to \infty$ and having only simple poles. Obviously, such matrix can be constructed
by using the arbitrariness of the rational functions~$\beta_n(k)$. Let the poles of $B(k)$
be located at the points $k = \rho_l$, $l = 1\dots d$.

Consider function
\begin{equation}
V(b,k) = U(b,k) B(k).
\label{eq0601}
\end{equation}
Note that
\begin{equation}
\frac{\ptl V(b, k)}{\ptl b} V^{-1} (b, k) = \frac{\ptl U(b, k)}{\ptl b} U^{-1} (b, k) \equiv S(b, k) .
\label{eq0601a}
\end{equation}
Thus, $V$ obeys ODE1 (\ref{eq0501}) with the same coefficient as $U$.

The key property of function $V$ is expressed by the following proposition.

\begin{proposition}
There exists function $R(b,k)$, which is rational as a function of $k$ for each $b$,
such that
\begin{equation}
V(b,k) = R(b,k) \, U(b, k).
\label{eq0602}
\end{equation}
\end{proposition}

Construct function $R$ as follows:
\begin{equation}
R(b, k) = V(b, k) \, U^{-1}(b,k) = U(b, k)\, B(k)\, U^{-1}(b,k)
\label{eq0603}
\end{equation}
Consider the behavior of $V(b,k)$ on the cuts $\Gamma_j(b)$.
Since $B$ commutes with all $H_j$,
\[
R(b, k^+) = U(b, k^-)\, H_j (k) \, B(k)\, H_j^{-1}(k) \, U^{-1}(b,k^-) =
\]
\[
U(b, k^-) \, B(k) \, U^{-1}(b,k^-) = R(b, k^-), \qquad k \in \Gamma_j (k).
\]
Thus, for each $b$ function $R(b,k)$ is a single--valued function of $k$. At infinity
$R(b,k) \to I$. Obviously, $R$ can only have  simple poles at $k = \rho_l$. Due to Liouville's
theorem, $R(b,k)$ should be a rational function of~$k$. Moreover, one can conclude
that $R$ has form
\begin{equation}
R(b, k) = I + \sum_{l = 1}^d \frac{r_l (b)}{k - \rho_l},
\label{eq0604}
\end{equation}
where $r_l (b)$ are some $N \times N$ matrix functions of $b$ defined in $\Omega$.

Construct the coefficient of ODE1 for $V$ using representation (\ref{eq0602}):
\begin{equation}
\frac{\ptl V}{\ptl b} V^{-1}  = R \,\frac{\ptl U}{\ptl b} \,U^{-1} \,R^{-1}
+ \frac{\ptl R}{\ptl b}\, R^{-1}
= R \,S \,R^{-1} + \frac{\ptl R}{\ptl b}\, R^{-1}.
\label{eq0605}
\end{equation}
Comparing (\ref{eq0605}) with (\ref{eq0601a}), conclude that
\begin{equation}
\frac{\ptl R(b, k)}{\ptl b} = S \, R - R\, S \equiv [S, R].
\label{eq0606}
\end{equation}
Equation (\ref{eq0606}) is the {\em global form\/} of ODE2. One can easily see that
this equation describes the evolution of $R$ but from the first glance it is
not clear how it can describe the evolution of $S$. However, we possess some additional information
about $R$ and $S$ (namely, both functions are rational with respect to $k$). This information is
enough to transform (\ref{eq0606}) into a {\em local form}, which is a closed set of
ordinary differential equations describing the evolution of $R$ and~$S$.

Substitute (\ref{eq0606}) and
\begin{equation}
S(b,k) = \sum_{j = 1}^p \frac{s_j (b)}{k - (k_j + b)}
\label{eq0607}
\end{equation}
into (\ref{eq0606}). Expand the right--hand side and left--hand side of (\ref{eq0606})
as a sum of simple fractions. Taking into account that
\[
\frac{1}{k - (k_j + b)}\frac{1}{k - \rho_l} =
\frac{1}{(k_j + b) - \rho_l}  \left(
\frac{1}{k - (k_j + b)}
-
\frac{1}{k - \rho_l}
\right)
\]
and considering the terms with each denominator separately, obtain equations
\begin{equation}
\frac{d r_l(b)}{d b} = \sum_{j=1}^p \frac{[ s_j(b), r_l(b) ] }{\rho_l - (k_j + b)},
\qquad
l = 1\dots d
\label{eq0608}
\end{equation}
and
\begin{equation}
 \sum_{l=1}^d \frac{[ s_j(b), r_l(b) ] }{\rho_l - (k_j + b)} = 0 .
\qquad
j = 1\dots p
\label{eq0609}
\end{equation}

System (\ref{eq0608}), (\ref{eq0609}) does not form a closed system of ordinary differential
equations for finding the unknown matrices $s_j(b)$, $r_j(b)$.
To make the system closed, consider (\ref{eq0609}) together with (\ref{eq0502a}). Formulate the problem
of finding of matrices $s_j$ provided that matrices $r_l$ are known. Equations (\ref{eq0502a}) provide
information about the eigenvalues of $s_j$, while (\ref{eq0609}) provide information
about the eigenvectors of $s_j$.
Namely, the eigenvalues of $s_j(b)$ are equal to the diagonal elements of
\[
\tilde F_j(b) = - \frac{1}{2 \pi i}  \log(F_j(k_j+ b)).
\]
where $F_j(k_j+ b)$ is a (known)
diagonal matrix composed of the eigenvalues of $H_j (k_j + b)$ (see (\ref{eq0302a})).
According to (\ref{eq0609}),
the eigenvectors of $s_j (b)$ coincide with the eigenvectors of the matrix
\[
R (b , k_j + b) = \sum_{l=1}^d \frac{ r_l(b) }{(k_j + b) - \rho_l}  .
\]

Define function ${\cal F}(X,Y)$ producing a matrix, whose eigenvalues coincide with the eigenvalues of
$X$, and the eigenvectors coincide with the eigenvectors of $Y$ (provided all eigenvalues of $Y$ are distinct).
The function $\cal F$ is defined ambiguously since the mapping between the eigenvalues of $X$ and eigenvectors
of $Y$ is not defined. I.e. $\cal F$ is defined up to a permutation of order $N$. If this ambiguity is
eliminated in a correct way,
\begin{equation}
s_j = {\cal F}(\tilde F_j(b) , R (b , k_j + b)).
\label{eq0610}
\end{equation}

Equations (\ref{eq0608}) together with (\ref{eq0610}) form a closed system of equations
for finding $r_l$ and~$s_j$. This system is non-linear. The system (\ref{eq0608}), (\ref{eq0610})
will be called the ODE2 (in the local form). Derivation of the ODE2 is the main result of this paper.
This result can be formulated in the form of the following theorem.

\begin{theorem}
Let there be a family of Riemann--Hilbert problems (\ref{eq0202a}) obeying the restrictions posed above, including the commutativity restrictions (\ref{eq0203}). The ODE1
for this family has the notation of (\ref{eq0501}).
Then there exist such matrices $r_l(b)$ and such a choice of the function $\cal F$ that matrices
$s_j(b)$, $r_l(b)$ obey the system (\ref{eq0608}), (\ref{eq0610}).
\end{theorem}

\subsection{Initial conditions for ODE2 and choice of function~${\cal F}$}

To make a numerical solution of ODE2 possible one should define the initial conditions and
eliminate the ambiguity of defining the function~$\cal F$. Since $U(b, k) \to I$
as $b \to i\infty$, one can conclude that
\begin{equation}
R(i \infty, k) = B(k),
\label{eq0611}
\end{equation}
where of course
\[
R(i \infty, k) \equiv \lim_{b \to i \infty} R(b, k).
\]
Thus, if
\begin{equation}
B(k) = I + \sum_{l= 1}^d \frac{t_l}{b- \rho_l}
\label{eq0611a}
\end{equation}
for some matrices $t_l$ (which are assumed to be known) then
\begin{equation}
r_l(i\infty) = t_l.
\label{eq0612}
\end{equation}
These relations play the role of initial conditions for the~ODE2.

To eliminate the ambiguity of definition of function ${\cal F}$,
we need to establish a correspondence between the eigenvectors of the
matrix $R(b, k_j + b)$ and the diagonal element of the (diagonal) matrix
$\tilde F_j (b)$. Again, consider large values of ${\rm Im}[b]$. For large
imaginary $b$ the values of the coefficients $H_j(k)$ approximately commute
with the common factor~$B(b)$. Therefore, the solution $U(b, k)$ near the
points $k = k_j + b$ approximately commutes with $B(b)$ or (which is the sam in asymptotic sense)
with $B(k_j + b)$. Thus, according to (\ref{eq0603}),
\[
R(b, k_j + b) \approx B(k_j + b).
\]
Using this relation,
one can establish a natural correspondence between the eigenvectors
of $B(k_j + b)$ and $R(b, k_j + b)$.
Then, the eigenvectors of $B(k_j + b)$ are by construction the eigenvectors of $H_j (k_j + b)$.
Thus, it is possible to establish a natural correspondence between the eigenvectors of $B(k_j + b)$
and $H_j (k_j + b)$. Finally, this gives correspondence between the diagonal elements of $F_j (k_j + b)$
and the eigenvectors of $R(b, k_j + b)$.

Thus, function ${\cal F}$ can be defined without ambiguity for large ${\rm Im}[b]$ and for other $b$
it can be defined by continuity.

\subsection{Invariance of ODE2 with respect to the choice of $B(k)$}

The choice of the factor $B(k)$ is not unique. Namely, if $B(k)$ obeys all restrictions
then a combination
\begin{equation}
B'(k) = \sum_{m = 0}^{N-1} g_m(k) B^m(k)
\label{eq2001}
\end{equation}
with rational scalar functions $g_m (k)$ also can be used as $B$, provided that $B'(k) \to I$
as $|k| \to \infty$ and $B'$ has only simple poles.
The form of ODE2 changes when $B$ is substituted by $B'$. Let us show that this substitution
does not change the solution $s_j (b)$. For this we remind that the system (\ref{eq0608}),
(\ref{eq0609}) is equivalent to (\ref{eq0606}). The invariance of $s_j$ is established by
the following proposition.

\begin{proposition}
Let $R$ be defined by (\ref{eq0603}), and
\begin{equation}
R'(b, k) = U(b, k)\, B'(k)\, U^{-1}(b,k),
\label{eq2002}
\end{equation}
where $B'$ is defined by (\ref{eq2001}). Let (\ref{eq0606}) be valid for some matrix $S$. Then
\begin{equation}
\frac{\ptl R'(b,k)}{\ptl b} = [S, R'].
\label{eq2003}
\end{equation}
\end{proposition}

First, note that it follows from (\ref{eq2002}) that
\begin{equation}
R'(k) = \sum_{m = 0}^{N-1} g_m(k) R^m(k)
\label{2004}
\end{equation}
Due to formal linearity of (\ref{eq2003}), it is sufficient to prove that
\begin{equation}
\frac{\ptl R^m(b,k)}{\ptl b} = [S, R^m].
\label{eq2003c}
\end{equation}
This can be easily proved by induction.


\section{Examples}

\subsection{Description of the numerical procedure}

The numerical procedure straightforwardly follows from Theorem~1 and Theorem~2.
Assume that matrices $H_j (k)$ are known explicitly, and let the matrix $B(k)$ be constructed and
represented in the form (\ref{eq0611a}).

First, ODE2 is solved
along the positive imaginary axis of $b$ from $i \infty$ to~$0$. In practice, ODE2 is solved not from
$b = i \infty$, but from $b = i L$, where $L$ is a large number playing the role of infinity.
At the ``infinite'' point $b = i L$ initial condition for ODE2 are set in the form of
\begin{equation}
r_l(i L) = t_l.
\label{eq0701}
\end{equation}
At the point $b = i L$ function ${\cal F}$ is constructed without ambiguity as follows.
According to the argument above and according to (\ref{eq0701}), for numerical solution
\begin{equation}
R(i L , k_j + iL) = B(k_j + i L).
\label{eq0702}
\end{equation}
Represent $B(k_j + iL )$ in the form (\ref{eq0303}), i.e.
\begin{equation}
B(k_j + iL) = P_* D_* P_*^{-1},
\label{eq0703}
\end{equation}
where $D_*$ is a diagonal matrix. Function $H_j(k_j + i L)$ can be represented in the form (\ref{eq0302a}),
i.e.
\[
P_*^{-1} H_j(k_j + i b) \, P_*
\]
should be a diagonal matrix. According to (\ref{eq0610}),
\begin{equation}
s_j(iL) = -\frac{1}{2\pi i} P_* \log( P_*^{-1} H_j(k_j + i b) \, P_* ) P_*^{-1},
\label{eq0704}
\end{equation}
where the branch of logarithm close to zero is taken. This procedure defines $s_j (i L)$
in a unique way.

Then ODE2, i.e.\ the system (\ref{eq0608}), (\ref{eq0610}) is solved numerically,
say by Runge--Kutta method, from $b = iL$ to $b = 0$.
On each step function ${\cal F}$ is chosen such that new values $s_j (b - i \delta )$ are close
to old values $s_j (b)$, i.e.\ such that $s_j (b)$ are continuous. As the result of this procedure,
the matrices $s_j$ are found at points covering the segment $(iL , 0)$ densely enough.

Next, ODE1 is solved to find $U(b,k)$. A set of points $k = z_n$ at which the solution
$U(k)$ will be found is selected. The initial conditions have form
\begin{equation}
U(iL , z_n) = I.
\label{eq0705}
\end{equation}
Equation (\ref{eq0502}) is solved from $b = i L$ to $b = 0$
for the values $U(b, z_n)$
along the imaginary axis of $b$ say
by Runge--Kutta method. As the result, the solution $U(k) = U(0, k)$ becomes known at the points
$k= z_n$.

One can see that the numerical procedure is rather simple. It consists of two solutions of ordinary differential
equations. If the segment $(iL, 0)$ is split into $N_b$ steps, and if there are $N_k$ points in the set $z_n$, then the first step takes $\sim N_b$ operations, and the second step takes $\sim N_k N_p$ operations.
This makes difference with results of \cite{arxiv2} where the first step takes $\sim N_p^2$ operations.

\subsection{Khrapkov's case}

It is important to show that the proposed technique is equivalent to the known method in the
simplest commutative case, namely in Khrapkov's case \cite{Khrapkov}.
Consider as an example a family of Riemann--Hilbert problems set on
$\Gamma_1 (b) = (k_1 + b , k_1 + i \infty)$ with the coefficient
\begin{equation}
H_1 (k) = g_0(k) I + g_1(k) \Lambda(k) ,
\label{eq0801}
\end{equation}
where
\[
\Lambda(k) = \left( \begin{array}{cc}
1 & k \\
k & -1
\end{array} \right),
\]
$g_0$ and $g_1$ are some algebraic functions such that $g_0 (k) \to 1$ as $|k| \to \infty$,
$g_1$ tends to zero as $|k|\to \infty$ not slower than $1/ |k|^2$.

A traditional solution of this problem is as follows. First, a solution $\tilde U (b,k)$ is constructed by
the formula \cite{Khrapkov}
\begin{equation}
\tilde U(b, k) = \exp (\bar \xi) \left(
\cosh \left( \sqrt{\phi(k)} \bar \eta \right) I +
\sinh \left( \sqrt{\phi(k)} \bar \eta \right) \frac{\Lambda(k)}{\sqrt{\phi(k)}}
\right),
\label{eq0802}
\end{equation}
\[
\phi(k) = k^2 +1 ,
\]
\begin{equation}
\bar \xi(b, k) = - \int \limits^{k_1 + i \infty}_{k_1 + b} \frac{\xi(\tau)}{k - \tau} d \tau,
\qquad
\bar \eta(b, k) = - \int \limits^{k_1 + i \infty}_{k_1 + b} \frac{\eta(\tau)}{k - \tau} d \tau,
\label{eq0803}
\end{equation}
\begin{equation}
\xi(k) = -\frac{1}{4 \pi i} \log \left( g_0^2(k) - \phi(k) g_1^2 (k) \right),
\label{eq0804}
\end{equation}
\begin{equation}
\eta(k) = -\frac{1}{4 \pi i \sqrt{\phi(k)}} \log \left(
\frac{g_0(k) + g_1(k) \sqrt{\phi(k)}}{g_0(k) - g_1(k) \sqrt{\phi(k)}}
 \right),
\label{eq0805}
\end{equation}

Solution $\tilde U$ obeys all conditions except the condition $U \to I$ at infinity. Instead,
for a fixed $b$
\begin{equation}
\tilde U(b , k) \to
\cosh( \zeta(b)) I + \sinh(\zeta(b)) \left(   \begin{array}{cc}
0 & 1 \\
1 & 0
\end{array} \right) \equiv Q(b) \qquad \mbox{as } |k| \to \infty,
\label{eq0806}
\end{equation}
\begin{equation}
\zeta(b)= - \int \limits^{k_1 + i \infty}_{k_1 + b}  \eta(\tau) d \tau.
\label{eq0807}
\end{equation}
Thus, one has to ``correct'' the behavior of $\tilde U$ by a left multiplication:
\begin{equation}
U(b , k) = Q^{-1} (b) \tilde U(b , k).
\label{eq0808}
\end{equation}

Let us consider the same problem from the point of view of the proposed method. One can check directly
\cite{arxiv2} that the auxiliary solution obeys ODE1 in a slightly modified form:
\begin{equation}
\frac{\ptl \tilde U(b, k)}{\ptl b} = \tilde S(b, k) \tilde U(b, k) ,
\qquad
\tilde S(b, k) =
\frac{\xi(b + k_1)}{k - (b + k_1)} I +
\frac{\eta(b + k_1)}{k - (b + k_1)} \Lambda(k)
.
\label{eq0809}
\end{equation}
Similarly, it can be checked that $Q (b)$ obeys equation
\begin{equation}
\frac{d Q}{d b} = \eta(b + k_1) \left(  \begin{array}{cc}
0 & 1 \\
1 & 0
\end{array} \right)\, Q(b).
\label{eq0809bb}
\end{equation}

Construct ODE1 for $U$. According to (\ref{eq0808}) the coefficient of this equation is equal to
\begin{equation}
S(b,k) = \frac{\ptl \tilde U}{\ptl b} = Q^{-1} \left(
\tilde S - \frac{d Q}{d b} Q^{-1}
\right) Q =
Q^{-1} \left(
\tilde S -
\eta(b + k_1) \left(  \begin{array}{cc}
0 & 1 \\
1 & 0
\end{array} \right)
\right) Q.
\label{eq0810}
\end{equation}
One can see that the coefficient has form of (\ref{eq0403}), i.e. for a fixed $b$ it is a rational
function of $k$ having a simple pole at $k = b+k_1$ and decaying at infinity.

Now consider ODE2. Select a function $B(k)$ commuting with both branches of $H_1(k)$, having only
simple poles and tending to $I$ at infinity.
For example one can choose
\begin{equation}
B(k) = I + \frac{1}{k^2 - 1} \Lambda(k)
\label{eq0811}
\end{equation}
with simple poles at $k = \pm 1$.
Define $\tilde R  = \tilde U B \tilde U^{-1}$. One can see that $B$ commutes with $\tilde U$,
and thus $\tilde R(b, k) = B(k)$.
Note that
\begin{equation}
[\tilde R , \tilde S] =0, \qquad \frac{\ptl \tilde R}{\ptl b} = 0.
\label{eq0812}
\end{equation}

Define $R$ as (\ref{eq0603}). It can be expressed as
\begin{equation}
R(b,k) = Q^{-1}(b) \tilde R(b, k) Q(b) = Q^{-1}(b) B(k) Q(b)
\label{eq0813}
\end{equation}
Taking into account (\ref{eq0812}) and (\ref{eq0810}) it is easy to show that
equation (\ref{eq0606}) i.e.\ ODE2 in the global form is valid for Khrapkov's matrix.

Let us write down ODE2 in the local form (for demonstration purposes).
Represent $\Lambda (k)$ in the form
\begin{equation}
\Lambda (k) = \sqrt{k^2 + 1} \, P(k) \left( \begin{array}{cc}
1 & 0 \\
0 & -1
\end{array}  \right) P^{-1} (k),
\label{eq0813a}
\end{equation}
\begin{equation}
P(k) = \left( \begin{array}{cc}
1 & 1 \\
\frac{\sqrt{k^2 +1}-1}{k } & \frac{-\sqrt{k^2 +1}-1}{k }
\end{array} \right).
\label{eq0814}
\end{equation}
Similarly,
\[
H_1 (k) = P(k) F_1 (k) P^{-1} (k),
\]
\begin{equation}
F_1 (k) =
\left( \begin{array}{cc}
g_0(k) + \sqrt{k^2 + 1} \, g_1 (k) & 0 \\
0 & g_0(k) - \sqrt{k^2 + 1} \, g_1 (k)
\end{array} \right)
\label{eq0815}
\end{equation}
Since $R$ has two poles ($\rho_1 = 1$, $\rho_2 = -1$), we need a system of equations describing
evolution of three matrices: $r_1 (b)$, $r_2 (b)$, and $s_1 (b)$. According to (\ref{eq0608}),
first two equations have form
\begin{equation}
\frac{d r_1 (b)}{db} = \frac{[s_1 (b) , r_1 (b)]}{1 - (k_1 + b)},
\qquad
\frac{d r_2 (b)}{db} = \frac{[s_1 (b) , r_2 (b)]}{-1 - (k_1 + b)},
\label{eq0816}
\end{equation}
The third equation has form of (\ref{eq0610}):
\begin{equation}
s_1 (b) = - \frac{1}{2 \pi i} {\cal F} \left(
\log(F_1 (b + k_1)) , \frac{r_1 (b)}{k_1 + b -1} + \frac{r_2 (b)}{k_1 + b + 1})
\right).
\label{eq0817}
\end{equation}

Initial conditions for (\ref{eq0816}) should be taken in the form (\ref{eq0612}). For this,
matrix $B$ should be represented as a sum of simple fractions. As the result, we get
\begin{equation}
r_1 (i \infty) = \frac{1}{2} \left( \begin{array}{cc}
1 & 1 \\
1 & -1
\end{array}\right) ,
\qquad
r_2 (i \infty) = \frac{1}{2} \left( \begin{array}{cc}
-1 & 1 \\
1 & 1
\end{array}\right)
\label{eq0818}
\end{equation}

Function ${\cal F}(X,Y)$ is implemented as follows. Let $X$ be a diagonal matrix. Matrix $Y$ is represented in the form
$Y = Y_1 \, Y_2 \, Y_1^{-1}$ numerically or analytically
($Y_2$ should be a diagonal matrix). The result is formed as
\begin{equation}
{\cal F}(X,Y) = Y_1 \, X \, Y_1^{-1} ,
\label{eq0819}
\end{equation}
or
\begin{equation}
{\cal F}(X,Y) = Y_1 \, X' \, Y_2^{-1}
\label{eq0820}
\end{equation}
where $X'$ is a matrix, whose diagonal elements are interchanged.
The choice between these two forms is made by the following rule.
For the point $b = i L$ where conditions (\ref{eq0818}) are set
matrix $Y$ has two eigenvectors, one of which is close
to
\[
a \left( \begin{array}{c} 1 \\ 1\end{array} \right) ,
\]
 and another one is close to
\[
a \left( \begin{array}{c} 1 \\ -1\end{array} \right)
\]
(see (\ref{eq0814})).
These eigenvectors are columns of $Y$. If the first column
corresponds to the vector of the first type, then
form (\ref{eq0819}) is chosen at this point. Otherwise,
form (\ref{eq0820}) should be chosen. At each new step function ${\cal F}$
is chosen to be approximately continuous.

\subsection{Factorization of Antipov's matrix}

Here we consider a more sophisticated (but also commutative) case previously
addressed in \cite{Antipov}. Matrix $G(k)$ is as follows:
\begin{equation}
G(k) = g_0 (k) I + g_1 (k) \Lambda(k),
\label{eq0901}
\end{equation}
where
\begin{equation}
\Lambda(k) = \left( \begin{array}{cc}
k^4 - \mu^4 & \alpha \mu^4 /\tau \\
\alpha \mu^4 /\tau & -k^4 + \mu^4
\end{array} \right)
\label{eq0902}
\end{equation}
\begin{equation}
g_0 (k) = \frac{(\psi(k)-\tau)(k^4 - \mu^4)-\alpha \mu^4}{\psi(k)(k^4 - \mu^4)},
\qquad
g_1 (k) = \frac{\tau}{\psi(k)(k^4 - \mu^4)},
\label{eq0903}
\end{equation}
\begin{equation}
\psi(k) = \sqrt{k^2 - (1 + 0i)^2}
\label{eq0904}
\end{equation}
$\mu$, $\tau$, $\alpha$ are some scalar constant physical parameters.
Notation (\ref{eq0904}) means that the only branch point in the upper half--plane is
$k_1 = 1$.

Matrix (\ref{eq0901}) is related to a problem of scattering by a screen composed of a rigid half--plane
and a flexible perforated sandwich half--plane. The boundary conditions for this problem were
derived in \cite{Leppington}. The problem was reduced to the Wiener--Hopf problem in \cite{Antipov}.

The problem belongs to the Khrapkov's class. The most important function for such problem is
\begin{equation}
\phi(k) =  \Lambda_{12} \Lambda_{21} -\Lambda_{11} \Lambda_{22} ,
\label{eq0905}
\end{equation}
having the property
\[
\Lambda^2(k) = \phi(k) I.
\]
In this case
\begin{equation}
\phi(k) = k^8 - 2 \mu^4 k^4 + \mu^8 (1 + \alpha^2 / \tau^2).
\label{eq0906}
\end{equation}
This function is a polynomial of degree~8. If a direct Khraphov's method \cite{Khrapkov} is applied
then a solution grows rapidly (faster than algebraically) at infinity. Therefore, the Moiseev's method
should be applied. The method has been outlined in \cite{Antipov}, however no numerical results have been
presented. An application of this method requires finding zeros of Riemann's theta function and Weierstrass'
kernel quadratures.

We apply the method developed above to this problem. The following values of parameters are taken for computations: \[
\mu = 2, \qquad \tau = 0.25, \qquad \alpha = 0.3.
\]
Apply the Hurd's method. There a single cut $\Gamma_1$ in the upper half-plane going from
$k_1 = 1$ to $1 + i \infty$. Denote by $g_0^+(k)$, $g_1^+(k)$ the values of $g_0(k)$, $g_1(k)$
on the right shore of the cut,
and by  $g_0^-(k)$, $g_1^-(k)$ the values on the left shore of the cut.
Note that these values are different due to the presence of the square root $\psi$.
The coefficient $H_1(k)$ describing the multiplicative jump on $\Gamma_1$
is equal to
\begin{equation}
H_1 (k) = G(k^-) G^{-1} (k^+) =
\frac{(g_0^- g_0^+ - \phi \, g_1^- g_1^+) I + (g_1^- g_0^+ - g_1^+ g_0^-) \Lambda}{(g_0^+)^2 - (g_1^+)^2 \phi }
\label{eq0907}
\end{equation}

Then, the function $B(k)$ is chosen. We can take it in the form
\begin{equation}
B(k) = I + \frac{1}{\xi(k)} \Lambda(k),
\label{eq0908}
\end{equation}
where $\xi(k)$ is a rational function. Since $\lambda(k)$ grows as $k^4$, we can take
$\xi (k)$ as a polynomial of 5th order, namely
\begin{equation}
\xi(k) = \prod_{l = 1}^5 (k - \rho_l).
\label{eq0909}
\end{equation}
The choice of $\rho_l$ can be done quite arbitrarily. We use the values
\[
\rho_1 =2 +  i, \quad
\rho_2 = 2 - i , \quad
\rho_3 =  - i ,  \quad
\rho_4 = -1+i , \quad
\rho_5 = -1 - i.
\]

\begin{figure}[ht]
\centerline{\epsfig{file=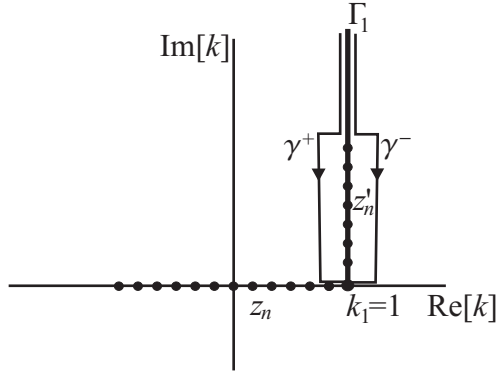}}
\caption{Contours $\gamma^+$ and $\gamma^-$}
\label{fig03}
\end{figure}

\begin{figure}[ht]
\centerline{\epsfig{file=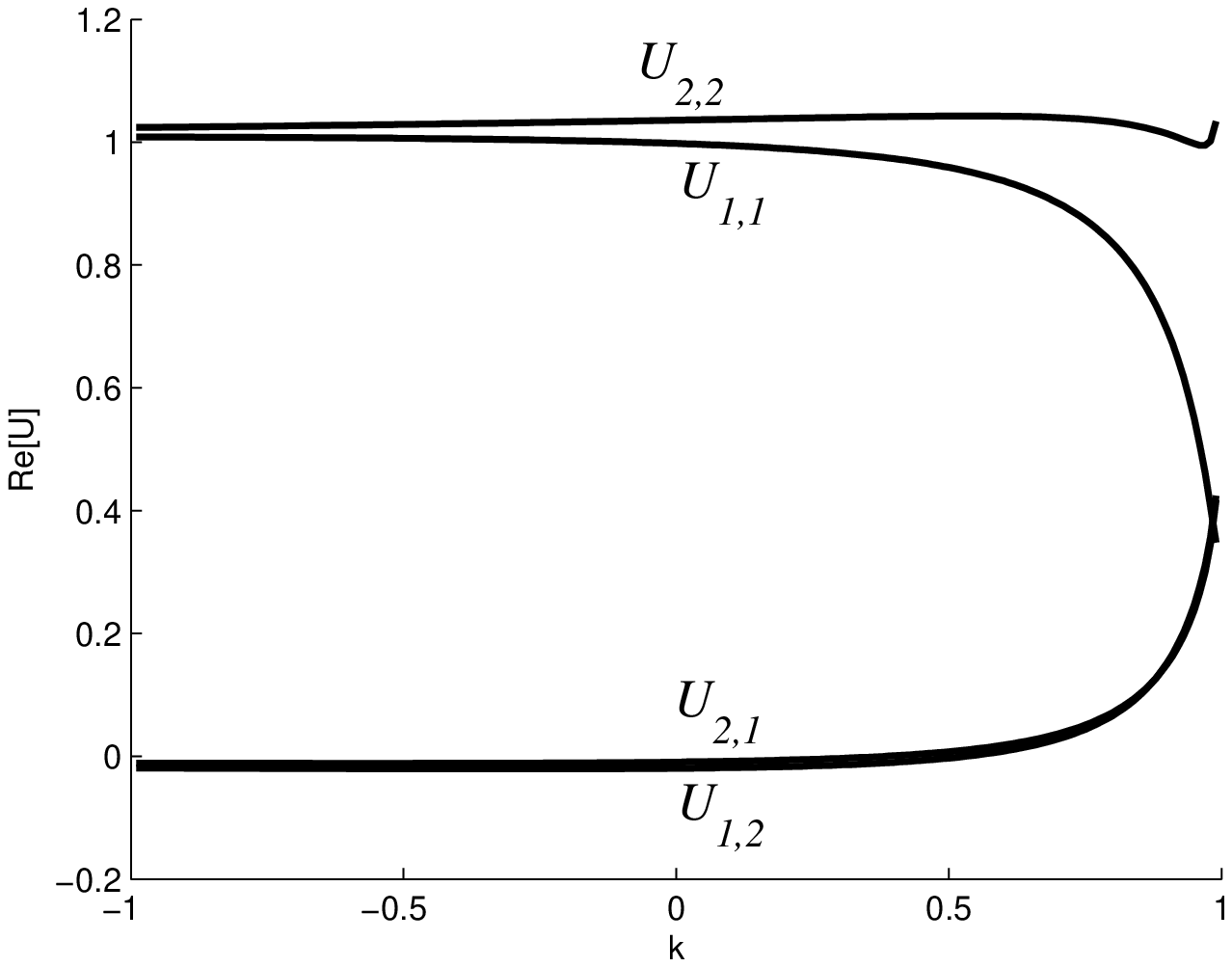,width=7cm}\epsfig{file=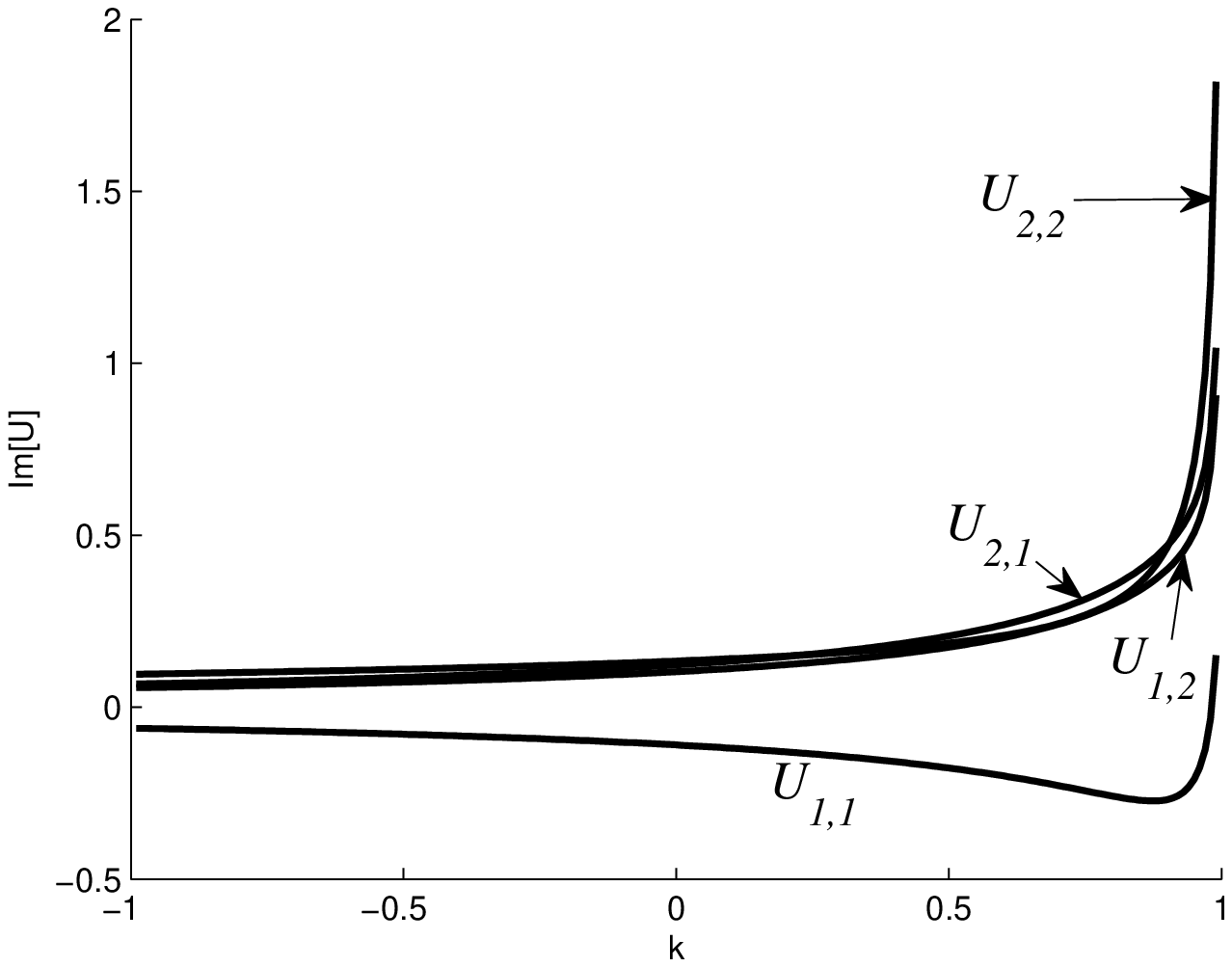,width=7cm}}
\caption{Solution $U(k)$ on the segment $k \in (-1,1)$} \label{fig04}
\end{figure}

The scheme outlined above is implemented. The set of the points of interest $k = z_n$
belong to the real segment $k \in (-1,1)$ (see Fig.~\ref{fig03}).
To determine them, it is necessary to find the values
$s_1 (b)$ for $b \in \Gamma_1$. These values are found by solving ODE2.
The result (i.e. the components of the matrix $U(0,z_n)$) is shown in Fig.~\ref{fig04}.

Besides finding the values $U(0, z_n)$ we perform a simple control of the whole procedure.
For this, we find the values $U(0,(z'_n)^+)$ and $U(0,(z'_n)^-)$ , where the points $z'_n$
belong to the cut $\Gamma_1$ (see Fig.~\ref{fig03}),
the values $U(0,(z'_n)^+)$ represent the right shore of the cut,
and the values $U(0,(z'_n)^-)$ represent the left shore of the cut. To determine the values on the shores we change the contour for solving ODE2 slightly. Namely, for the values $U(0,(z'_n)^+)$ we chose contour $\gamma^+$
in Fig.~\ref{fig03}, and for values $U(0,(z'_n)^-)$ we chose contour~$\gamma^-$. After that, we compute the
combination $U^{-1}(0,(z'_n)^-) U(0,(z'_n)^+) M^{-1}(z'_n)$. In the ideal case this matrix should be equal to
$I$, therefore its deviation from $I$ can be taken as a measure of relative accuracy of the computation. It has been found that the relative accuracy of the computation used to be of order~$10^{-4}$.

\section{Conclusion}

A new method for matrix factorization in the commutative (Moiseev's) case is developed. The method
is numerical, but it is based on two analytical properties of the factorization problem. It is applicable
to algebraic matrices having the property of branch--commutativity, i.e.\ the matrices, whose values corresponding
to different sheet over the same affix commute.

The factorization problem is transformed by Hurd's procedure into a Riemann--Hilbert problem on a set of cuts.
Then the Riemann--Hilbert problem is embedded into a family of Riemann--Hilbert problems indexed by a variable~$b$. The solution as a function of $b$ is described by ordinary differential equation (ODE1) with an unknown coefficient~$S$. This coefficient is found by solving another ordinary differential equation, ODE2. Initial
conditions for ODE1 and ODE2 are formulated. It is shown that the proposed procedure in the Khrapkov's case is equivalent to the standard solution. Moreover, it is shown that the new procedure is applicable to Antipov's
matrix, and it does not lead to Jacobi's inversion problem, which is not easy to implement.

In more general (non-commutative) cases ODE2 should be replaced by an OE-equation described in \cite{arxiv2}.

\section*{Acknowledgements}

The work is supported by RF Government grant 11.G34.31.0066,
``Scientific school'' grant 2631.2012.2, RFBR grant 12-02-00114.
Author is grateful to Prof.\ R.V.\ Craster for valuable help.


\begin{thebibliography}{99}

\bibitem{Noble}
B. Noble, {\em Methods based on the Wiener--Hopf technique\/} (Pergamon Press, London 1958).

\bibitem{Moiseev}
 N. G. Moiseyev, Factorization of matrix functions of special form,
{\em Soviet Math. Dokl.} {\bf 39}  (1989) 264--267. 

\bibitem{Zverovich}
 E. I. Zverovich, Boundary value problems in the theory of analytic
functions in Holder classes on Riemann surfaces. {\em Russian Math.
Surveys\/} {\bf 26} (1971) 117--192 . 

\bibitem{Antipov}
Y. A. Antipov and V. V. Silvestrov, Factorization on a Riemann surface in scattering theory,
{\em Quart. J. Mech Appl. Math.\/} {\bf 55} (2002) 607--654.

\bibitem{Khrapkov}
A. A. Khrapkov, Certain cases of the elastic equilibrium of an infinite wedge with a nonsymmetric
notch at the vertex, subjected to concentrated forces, {\em J. Appl. Math. Mech. (PMM)\/} {\bf 35} (1971) 625–-637. 



\bibitem{Jones}
D. S. Jones, Commutative Wiener--Hopf factorization of a matrix, {\em Proc. R. Soc.~A\/}
{\bf 393} (1984) 185-–192.

\bibitem{Daniele}
V. G. Daniele, On the solution of vector Wiener--Hopf equations
occuring in scatering problems. {\em Radio Science\/} {\bf 19}  (1984)
1173--1178. 

\bibitem{AbrahamsCom}
B. H. Veitch and
I. D. Abrahams, On the commutative factorization
of $n\times n$ matrix Wiener-–Hopf kernels with distinct eigenvalues. {\em Proc. R. Soc.~A\/} {}
{\bf 463} (2007) 613--639.

\bibitem{Gakhov}
F. D. Gakhov, Riemann's boundary problem for a system of n pairs of functions, {\em Usp.\ Math.\ Nauk\/} {\bf 7} (1952) 3–-54.

\bibitem{Abrahams}
 I. D. Abrahams, On the solution of Wiener-Hopf
problems involving noncommutative matrix
kernel decompositions, {\em SIAM J. Appl. Math.} {\bf 57} (1997) 541--567.

\bibitem{Hurd}
R. A. Hurd,  the Wiener--Hopf--Hilbert method for diffraction problems, {\em Can. J. Phys.} {\bf 54} (1976) 775--780.

\bibitem{Chebotarev}
G. N. Chebotarev, On closed-form solution of a Riemann boundary value problem for $n$ pairs of functions, {\em Uchen. Zap. Kazan. Univ.} {\bf 116} (1956) 31–-58.

\bibitem{arxiv1}
A. V. Shanin, E. A. Doubravsky, Criteria for commutative factorization of a class of algebraic matrices,
arXiv:1211.4424.

\bibitem{Lancaster}
P. Lancaster. {\em Theory of matrices} (Academic Press, New-York --- London, 1969).

\bibitem{arxiv2}
A. V. Shanin, An ODE--based approach to some Riemann--Hilbert problems motivated by wave diffraction,
arXiv:1210.1964

\bibitem{Leppington}
F. G. Leppington, The effective boundary conditions for a perforated sandwich panel in a compressible
fluid, {\em Proc. R. Soc. Lond. A}{\bf 427} (1990) 385--399.


\end{thebibliography}
\end{document}